\documentclass[dvipsnames,11pt,reqno]{amsart}

  \usepackage{amssymb,mathtools,amsthm}
  \usepackage[a4paper,left=3cm,right=3cm,top=3cm,bottom=3cm]{geometry}
  \usepackage[T1]{fontenc}
  \usepackage[utf8]{inputenc}
\usepackage{microtype}
\usepackage{enumitem}

\newcounter{n}
\numberwithin{n}{section}
\theoremstyle{plain}
  \newtheorem{lemma}[n]{Lemma}
  \newtheorem{theorem}{Theorem}
  \newtheorem{proposition}[n]{Proposition}
  
\theoremstyle{definition}
  
  \newtheorem*{definition*}{Definition}
  \newtheorem{remark}[n]{Remark}

  \usepackage{tikz}
  \usepackage{subcaption}

\usepackage{hyperref}
\definecolor{colorlinks}{RGB}{0, 24, 168}
\definecolor{colorcites}{RGB}{124, 10, 2}
\hypersetup{
    colorlinks=true,
    linkcolor=colorlinks,
    citecolor=colorcites,
    urlcolor=colorlinks,
    pdfborder={0 0 0}
}

\renewcommand\phi\varphi
\renewcommand\epsilon\varepsilon
\DeclareMathSymbol{\shortminus}{\mathbin}{AMSa}{"39}

\begin{document}

 \title[Non-reversible stationary measures on trees]{Non-reversible stationary states for majority voter and Ising  dynamics on trees}

 \makeatletter
 \@namedef{subjclassname@2020}{\textup{2020} Mathematics Subject Classification}
 \makeatother

 \author{Piet Lammers}
 \address{Institut des Hautes \'Etudes Scientifiques}
 \email{lammers@ihes.fr}
 \author{Fabio Toninelli}
 \address{Technische Universit\"at Wien}
 \email{fabio.toninelli@tuwien.ac.at}
 
\begin{abstract}
We study three Markov processes on infinite, unrooted, regular trees: the stochastic Ising model (also known as the  Glauber heat
bath dynamics of the Ising model), a majority voter dynamic, and a coalescing
particle model. In each of the three cases the tree exhibits a preferred
direction encoded into the model.  For all three models, our main
result is the existence of a stationary but non-reversible measure. For the
Ising model, this requires imposing that the inverse temperature is large and
choosing suitable non-uniform couplings, and our theorem implies the existence
of a stationary measure which looks nothing like a low-temperature
Gibbs measure.  The interesting aspect of our results lies in the fact that the
analogous processes do not have non-Gibbsian stationary measures on $\mathbb
Z^d$, owing to the amenability of that graph. In fact, no example of a
stochastic Ising model with a non-reversible stationary state was known to date.

\end{abstract}

\maketitle

\setcounter{tocdepth}{1}

\section{Introduction}

\subsection{Interacting particles on trees}
We study three Markovian interacting particle processes: a model of
\emph{coalescing particles} which move at random times and join when they meet (Fig.~\ref{sub:particle}), a \emph{majority voter model} (Fig.~\ref{sub:voter}), and the
\emph{stochastic Ising model} with non-homogeneous coupling constants (see Fig.~\ref{sub:ising}; we also refer to the dynamics of this model as
the Ising dynamics, Glauber dynamics, or a variation thereof).  In all
three cases, the elementary degrees of freedom, or \emph{spins}, are placed on
the vertices of an infinite regular tree.  For the first two models this tree is
in fact a directed tree.
Our goal here is to exhibit some intriguing and somewhat surprising phenomena,
which \emph{do not occur} on finite-dimensional lattices and which we believe
root in the non-amenability and non-unimodularity of infinite directed trees
(see~\cite{MR3616205} for an introduction on these notions).

\subsection{The stochastic Ising model on Euclidean lattices}
\begin{figure}
    \centering
    \includegraphics{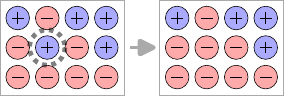}
    \caption{The Ising model of a portion of the two-dimensional square lattice,
    as a proxy for the Euclidean plane.
    Each vertex is endowed with a Poisson clock. When a clock rings, the spin
    is resampled according to the Ising model.
    Broadly speaking, this means that the spin is set to agree with the majority of its neighbours
    (breaking ties by flipping coins) most of the time,
    while occasionally the vertex chooses its spin independently of the other spins on the graph.
    In this case, the $+$ spin rings and turns into $-$,
    agreeing to the majority of its four neighbours after the clock rings.
    In 2D, it is known that the measures which are stationary for this dynamic coincide exactly
    with the set of Gibbs measures of the Ising model.}
    \label{fig:euclid}
\end{figure}

To illustrate our point, consider a locally finite graph $G=(V,E)$ and suppose
that the vertices in $V$ are endowed with i.i.d.\ Poisson clocks.  The
\emph{stochastic Ising model} (also known as the Glauber heat bath dynamic of the Ising model) is the Markov process on $\{+,-\}^V$ where each
\emph{spin} $\sigma_v$ with $v\in V$ is updated each time the clock at that
vertex rings, with the new spin value sampled from the specification of the Ising model
(see Fig.~\ref{fig:euclid} for an illustration and Section~\ref{sec:Ising} for details).
On
the one hand, Gibbs measures for the Ising model are well-known to be stationary and
reversible for this dynamic \cite{Liggett1}.  On the other hand, it is known that stationary
measures for the stochastic Ising model are also Ising Gibbs measures, at least
when:
\begin{enumerate}
    \item The underlying graph is the square lattice $\mathbb Z^d$ in dimension
    $d=1,2$~\cite{holley1977one};
    \item The underlying graph is the square lattice $\mathbb Z^d$ in dimension
    $d\geq 3$, and we restrict to measures which are invariant under the
    symmetries of the graph; see for example~\cite[Chapter~IV.5]{Liggett1}.
\end{enumerate}
In general, dynamics of these types have been widely studied on Euclidean
lattices in dimension $d\geq 2$ such as the square lattice, see for instance the
classical books \cite{Liggett2,Liggett1}.  To the best of our knowledge, it is
not known in dimension $d\geq 3$ if there exist stationary measures for the
dynamic which are not Gibbs measures for the Ising model.  It is, however,
expected that the situation is no different for $d\geq 3$ than for $d=1,2$.
Thus, a natural question is to ask: \emph{are there graphs for which there
exists a stationary measure which is non-reversible for the stochastic Ising model, and therefore not Gibbs?}

\subsection{The failure of entropy arguments on the tree}
\begin{figure}
    \begin{subfigure}{\textwidth}
        \centering
    \includegraphics{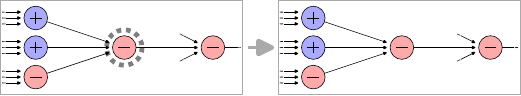}
    \subcaption{The Ising model}
    \label{sub:ising}
    \end{subfigure}
    \begin{subfigure}{\textwidth}
        \centering
    \includegraphics{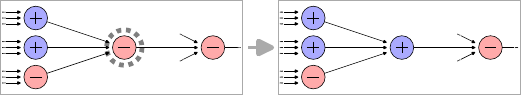}
    \subcaption{The voter model}
    \label{sub:voter}
    \end{subfigure}
    \begin{subfigure}{\textwidth}
        \centering
    \includegraphics{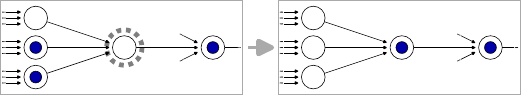}
    \subcaption{The particle model}
    \label{sub:particle}
    \end{subfigure}
    \caption{The three models under consideration. For the 
    Ising model, the illustration shows the unlikely event
    that the updated spin goes against the majority of its three children---the
     coupling constants are chosen such that the spins of the children weigh
     much heavier than the spin of the parent.}
    \label{fig:allmodels}
\end{figure}

On the square lattice (see also~\cite{higuchi1975some,moulin1977free}), the
results are, informally speaking, proved as follows (see Remark \ref{rem:Zd} below for more details on a similar argument): for a measure $\mu$ to be
stationary, its entropy has to be constant in time. On the other hand, the time
derivative of the entropy restricted to a large box is the sum of a \emph{bulk
term} proportional to the box volume, which is non-zero under the assumption
that $\mu$ is not reversible, plus a \emph{boundary term}. In finite dimension,
the boundary term cannot compensate for the bulk term because it is of smaller
order, thus contradicting non-reversibility. The same type of argument can
obviously not work on a non-amenable graph, but it is not a priori clear whether
this is a (technical) shortcoming of the proof or reflects a phenomenon of genuine interest. Our results
indicate that the latter is true: in particular, in Theorem~\ref{thm:ising}
below we show that the stochastic Ising model on a regular infinite tree admits
at least one non-reversible stationary measure $\mu$, if the inverse temperature
is large and for a suitable, space-dependent, choice of coupling constants (Fig.~\ref{sub:ising}).  We emphasize
that $\mu$ does not at all resemble a low-temperature Ising Gibbs measure: on
certain \emph{layers} (or \emph{generations}) of the tree the spins are
approximately i.i.d.\ and mix exponentially fast in time for the Ising dynamics.

An even neater phenomenon is exhibited by the majority voter model described in 
Subsection~\ref{sec:voter}, where the transition rates are the same at each
vertex. In this process, each spin or \emph{opinion} is updated at times
determined by i.i.d.\ exponential clocks and the new value of the spin of a
vertex is that of the majority of its neighbours which have an edge that is
incoming to the vertex (Fig.~\ref{sub:voter}).  We study this model on a directed tree where each
vertex has one outgoing edge and $d\in 1+2\mathbb N$ incoming edges
($d$ is odd in order to avoid ties).  For the analogous process on a
finite-dimensional lattice of odd degree, it is easy to see that any stationary,
shift-invariant measure is supported on completely frozen configurations (that is, absorbing states) such
that each spin agrees with the majority of its neighbours. In particular, such
stationary measures are (trivially) reversible. In contrast, for the majority
voter model on the tree we claim in Theorem~\ref{thm:nonfrozenvoters} that there
exists an automorphism-invariant stationary measure $\mu$ that is \emph{non-reversible}
and \emph{very far from frozen}: vertices in the same generation of the tree
have i.i.d.\ opinions, and the opinion of a single vertex mixes exponentially
fast in time.

Let us finally consider the above mentioned model of coalescing particles on
the tree (Fig.~\ref{sub:particle}).
In this model, each vertex is either empty or it has a particle on it.
When a clock at some vertex $x$ rings, all particles at children of $x$
are moved to $x$.
Particles merge when they meet.
For this model, we prove similarly (see Theorem~\ref{thm:info}) the existence of a
nontrivial  shift-invariant stationary measure $\mu$, while for the
analogous process on $\mathbb Z^d$ the only stationary measure would be the
Dirac mass on the completely empty configuration.

The stationary measures $\mu$ of the three models mentioned above, being non Gibbsian, do not have
an explicit description as the exponential of an energy function; rather, they are defined via a
limiting recursive procedure that is made possible by the directed tree
structure of $G$.

\subsection{Further questions}
Our results raise several rather intriguing questions regarding the Ising model.
One concerns the basin of
attraction of the non-reversible stationary measures: which initial
distributions evolve asymptotically into non-reversible stationary measures,
and which evolve into reversible Gibbs measures?  Another question is to ask whether the stochastic Ising
model on the regular infinite $d$-ary tree with constant couplings and large
enough $\beta$, or the majority voter model on the \emph{undirected} regular
$d$-ary tree (say with $d$ odd, to avoid ties), admit non-reversible stationary
states.

\subsection{Related literature}
While we are not aware of previous results on the existence of non-reversible
stationary measures for Ising Glauber dynamics, there are several works in the
literature that are related to our work.  The Glauber dynamics of the Ising
model on trees has widely been studied, see for
instance~\cite{berger2005glauber,martinelli2007fast,martinelli2004glauber}.
However, the main focus of these references is on finite (large) trees, where
the finite-volume Gibbs measure is the unique stationary one, and the main
results are estimates of the mixing or relaxation times of the dynamics. 
The relationship between reversibility, stationarity, and Gibbsianity
has been studied in several works, see
for example~\cite{kunsch1984non,kunsch1984time,maes1994non}.
We also refer to~\cite{durrettfinalref} for an analysis of the contact
process on the tree.
As
for the voter model (or zero-temperature Ising dynamics) on infinite, undirected
regular trees we mention for instance~\cite{howard2000zero,kanoria2011majority,eckner2014fixation}.
In the latter two works, one estimates the value of the threshold for consensus.
The \emph{threshold for consensus} is defined to be the
minimal bias in the initial i.i.d.\ opinion distribution that induces asymptotic
consensus for large times. The analogous question on $\mathbb Z^d$ is analysed
in~\cite{fontes2002stretched,morris2011zero}.

The majority voter model in our article is closely related to the \emph{Boolean
decision tree} with the \emph{recursive majority-of-three} algorithm implemented
at each vertex.  The difference is that in this model, there is a base layer
with initial votes, and the elections are held once and in order. In our model
there is no base layer, and elections are held at random times.  The Boolean
decision tree is studied in~\cite{jayram2003two} and~\cite{magniez2016improved},
where the authors find bounds on the computational complexity of these
functions.  The computational complexity is large because most of the opinions
at the base layer must be revealed in order to determine the final outcome of
the vote.  This is also closely related to the concept of \emph{noise
sensitivity}, see for example the book of Garban and
Steif~\cite{garban2014noise}. Without giving a full, rigorous explanation, the
authors of the current article believe that the existence of a non-frozen state 
for the majority voter model is closely related to the high computational
complexity and noise sensitivity of the Boolean decision tree with the recursive
majority-of-three algorithm.

Finally, let us mention that our techniques are heavily based on recursive formulae
which appear abundantly in the study of Gibbs measures on trees.
In particular, we recall the Zachary construction~\cite{zachary1983countable,zachary1985bounded},
which is a general construction for classifying Gibbs measures on trees
which exhibit the Markov property;
refer also to the book of Rozikov~\cite{rozikov2013gibbs}
and the more recent work of Shriver~\cite{shriver2020free}
on Gibbs measures on Cayley trees and 
the work of Sullivan~\cite{sullivan1975markov} and Vasilyev~\cite{vasilyev1978bernoulli} and the reference book of Georgii~\cite[Chapter~12]{georgii2011gibbs}
for general Gibbs measures.

\section{Definitions and main results}
\label{sec:def}

Trees and Poisson processes appear in all our models,
which is why we define them first.

\begin{definition*}
[Regular directed trees]
Let $\mathcal T=\mathcal T^d=(\mathbb V,\mathbb E)$ denote the directed unrooted tree in which each vertex has one parent and $d$ children.
Let $p(x)$ denote the parent of some vertex $x\in\mathbb V$,
let $e(x)$ denote the edge directed from 
$x$ to its parent,
and write $\mathcal C(x)$ for the set of children of $x$.
Also write $\mathcal D(x)$ for the set containing $x$ and all its descendants.
Choose some bi-infinite path $(x_k)_{k\in\mathbb Z}$ through $\mathcal T$ such that each $x_k$ is the parent of $x_{k-1}$.
There is a natural partition $(\mathbb V_k)_{k\in\mathbb Z}$ of $\mathbb V$ into \emph{layers},
which is defined such that $x_k\in\mathbb V_k$ and such that $p(x)\in\mathbb V_{k+1}$
whenever $x\in\mathbb V_k$.
Write $\operatorname{Aut}(\mathcal T)$
for the group of automorphisms 
of the directed tree $\mathcal T$.
Finally, fix some automorphism
 $\tau\in\operatorname{Aut}(\mathcal T)$ which has the property that each $x_k$
 is mapped to $x_{k-1}$.
\end{definition*}

\begin{definition*}[Poisson clocks]
Each vertex of the tree is endowed with a Poisson clock.
Formally, let $\mathbb P$ denote a probability measure in which 
$N$ is a Poisson random measure on $\mathbb V\times\mathbb R$
with intensity $\gamma\times\lambda$,
where $\gamma$ and $\lambda$
denote the counting measure and the Lebesgue measure on $\mathbb V$
and $\mathbb R$ respectively.
Denote the times at which the clock at the vertex $x$ rings 
by the sequence $(X^x_k)_{k\in\mathbb Z}$.
Formally,
this means that $(X^x_k)_{k\in\mathbb Z}$
is almost surely the unique increasing sequence of real numbers such that
\[
    N(x,\cdot)=\sum_k\delta_{X^x_k},
\]
and such that $X^x_1$ is the first positive number in the sequence.
\end{definition*}

Note that the above definition is valid for almost all $N$,
which is sufficient for our purposes.
The random measure $N$ satisfies a range of local finiteness conditions
which are natural and standard.
Such local finiteness conditions are important for the models below to be
well-defined.

\subsection{The coalescing particles model}
\label{subsec:main:info}

We term our first model the \emph{coalescing particles model} because it
can be viewed as a system of particles on
the tree with coalescence: particles coalesce when they happen to be
on the same vertex.  For the exposition of our results, we always
define this model on $\mathcal T^2$: the directed tree in which each
vertex has one parent and two children.  An identical analysis works
for a tree of higher degree, except that the exact formula in
Theorem~\ref{thm:info} below is harder to guess.  At each point in
time, each vertex has two possible states: either there is a
particle at that site, or not.  Informally, the dynamics
are defined such that whenever the clock at some vertex $x$ rings, all
particles at the children of $x$, are passed to $x$.
The particles merge into one if multiple particles are present. Formally, for
$\omega\in\{0,1\}^\mathbb V$, write $\omega^x\in\{0,1\}^\mathbb V$ for
the unique configuration such that
\begin{equation}
    \label{defcoadynmon}
    \omega^x_y=
    \begin{cases}
        0&\text{if $y\in\mathcal C(x)$,}\\
        1&\text{if $y=x$ and $\omega_z=1$ for some $z\in{\mathcal C(x)\cup\{x\}}$,}\\
        \omega_y&\text{otherwise.}
    \end{cases}
\end{equation}
When a Poisson clock at some site $x$ rings,
the state $\omega$ is replaced by the state $\omega^x$.
Observe that this update is local, which means that one can even make sense of 
this process on the infinite tree.
Write $L$ for the generator corresponding to this Markov process and $\omega(t)=\{\omega_x(t)\}_{x\in \mathbb V}$ for the configuration at time $t$ (analogous notations are used later for the majority voter and Ising dynamics).
Note also that the sample space $\{0,1\}^\mathbb V$ has a natural partial ordering $\leq$ where  $\omega\le \omega'$ if and only if $\omega_x\le \omega'_x$ for all $x\in\mathbb V$. We take the product topology on the sample space.

{
The dynamics in~\eqref{defcoadynmon} preserve the partial ordering: if $\omega\leq\tilde\omega$,
then $\omega^x\leq\tilde\omega^x$. Let $\mu_t$ denote the distribution on $\{0,1\}^\mathbb V$
after running the dynamics for time $t\in[0,\infty)$ started from the unique maximal configuration
where there is a particle at every vertex.
By the previous monotonicity, the map $t\mapsto \mu_t$ is stochastically decreasing with respect to $\leq$.
This means that it tends to a limit distribution $\mu_\infty$ on $\{0,1\}^\mathbb V$,
namely the unique maximal $L$-stationary distribution.
}

\begin{theorem}[Nontriviality of the coalescing particle process]
    \label{thm:info}
    There is a unique maximal probability distribution $\mu_\infty$ on $\{0,1\}^\mathbb V$
    which is stationary under $L$.
    This distribution is nontrivial in the sense that 
    it is not the Dirac measure on the outcome where no vertex has a particle.
    Moreover, it has the property that the probability
    of new particles arriving at a vertex $x$ in a time interval of length 
    $T$ is exactly equal to
    \[
        \rho_\infty(T):=1+6\frac{\beta e^{-T}}{(\beta e^{-T}-1)^2}
    ;\quad\beta=\sqrt 3 - 2<0.
    \]
\end{theorem}

\begin{remark}
    The nontrivial distribution in the previous theorem is not reversible
    under the dynamics,
    because the direction of the flow of particles respects the orientation of
    the edges.
    Moreover, this means that large portions of the tree do not communicate 
    with one another, which induces independence for such respective parts.
    More precisely, if two vertices $x$ and $y$ have the property that their \emph{parents}
    are incomparable with respect to the descendance structure on the tree (in the sense that no
    vertex is a descendant of the other),
    then the dynamics defined in~\eqref{defcoadynmon} preserve the independence between
    the descendants of $x$ and the descendants of $y$.
    In particular, this independence holds true for $\mu_\infty$ since it is defined by starting the
    dynamics
    from a deterministic configuration (which makes any two events independent).
    Consider, for example, a single layer of vertices.
    This layer has a natural partition, where each member contains precisely 
    one vertex and its siblings (that is, the vertices having the same parent).
    Then the distribution of $\mu_\infty$ is i.i.d.\ over this partition.
\end{remark}

\begin{remark}
  \label{rem:Zd}
    One can define a similar process on the undirected graph $\mathbb
    Z^d$, where initially there are particles at all vertices, and
    where at each Poisson clock ring at some vertex $x$, all particles
    at neighbours of $x$ are drawn towards $x$ and coalesce.
    Let $\mu_t$ denote the distribution on $\{0,1\}^{\mathbb Z^d}$
    after running the dynamics for time $t$.
    Exactly like before, $t\mapsto \mu_t$ is stochastically decreasing
    and thus tends to some limit $\mu_\infty$,
    which is the unique maximal stationary distribution for the dynamic.
    Note that the initial configuration and the dynamics are shift-invariant,
    which leads us to conclude that $\mu_\infty$ is also shift-invariant.

    Contrasting Theorem~\ref{thm:info},
    we claim that in fact $\mu_\infty$ equals a Dirac mass at the empty configuration
    (that is, all particles have vanished).
    Assume the claim to be false, which implies that the density
    $\delta$ (well-defined through shift-invariance) is strictly positive.
    We claim that the time derivative of $\delta$ is in fact strictly negative,
    contradicting stationarity.
    Notice that if $\delta>0$, then there is also a positive density $\delta'>0$
    of edges containing two particles.

    The time derivative of the
    average number of particles in a large cubic box $\Lambda$ equals
    $-2\delta'$ times the number of edges in the box 
    (accounting for the event that two particles merge),
    plus a boundary term of order $|\partial \Lambda|$
    (accounting for particles flying in or out of the box).
    If the box is large enough then the first term dominates,
    which implies strict negativity of the time derivative of $\delta$.

    Similar considerations apply to the directed
    version of the square lattice where all edges are oriented towards
    the positive direction of their corresponding coordinate axis.
    This means that Theorem~\ref{thm:info} relies on the tree
    structure in a nontrivial way.
\end{remark}

\subsection{The majority voter model}
\label{sec:voter}
We now come to our second model,
the \emph{majority voter model}.
At each point in time,
each vertex has a \emph{spin} or \emph{opinion},
which is $+$ or $-$.
Whenever a clock rings at some vertex,
its opinion is replaced by the majority of the opinions of its children.
To avoid tied votes, we define this model on $\mathcal T^3$,
the directed tree in which each vertex has one parent 
and three children.
Formally, for $\omega\in\{+,-\}^\mathbb V$, write $\omega^x\in\{+,-\}^\mathbb V$
for the unique configuration such that
\[
    \omega^x_y=\begin{cases}
        \omega_y&\text{if $y\neq x$,}\\
        +&\text{if $y=x$ and $\omega_z=+$ for at least two children $z$ of $x$,}\\
        -&\text{if $y=x$ and $\omega_z=-$ for at least two children $z$ of $x$.}
    \end{cases}    
\]
Thus, when a clock at some vertex $x$ rings,
the configuration $\omega$ is replaced by $\omega^x$.
Write $L$ for the corresponding generator.
A configuration is called \emph{frozen} if the opinion of each vertex
equals the majority of the opinions of its children.
A measure is called \emph{frozen}
whenever it is supported on frozen configurations.
Frozen measures are automatically stationary and reversible for the dynamic.

\begin{theorem}[Nontriviality of the majority voter model on the directed tree]
    \label{thm:nonfrozenvoters}
    There exists a stationary, $\operatorname{Aut}(\mathcal T)$-invariant Markov 
    process with generator $L$, which is supported 
    on non-frozen states, and which has the following two properties:
    \begin{enumerate}
        \item The opinion of each vertex $x\in\mathbb V$
        mixes exponentially fast in time, that is, 
        $\mathbb E[\omega_x(0)\omega_x(t)]:=\mathbb P[\{\omega_x(0)=\omega_x(t)\}]-\mathbb P[\{\omega_x(0)\neq\omega_x(t)\}]$
        decays exponentially fast in $t$, where $\mathbb E$ is the expectation with respect to the distribution of the stationary process.
        \item Its restriction to a single layer is
        i.i.d.\ over the vertices in that layer.
    \end{enumerate}
\end{theorem}

\begin{remark}
    This distribution is also not reversible.
    Indeed, if the process is run for a positive amount of time and then reversed,
    then the vertices change their opinion from majority to minority. 
    It is straightforward to generalise this theorem to directed trees
    $\mathcal T^d$ with $d$ an odd integer above three,
    at the cost of slightly complicating the calculations.
\end{remark}

\begin{remark}
In analogy with Remark~\ref{rem:Zd}, one can define the majority dynamics on a
$d$-dimensional periodic lattice of odd degree (to avoid ties), for instance by taking $\mathbb Z^d$ and removing one edge per vertex, in a periodic way. In this case, it
is easy to see that shift-invariant, stationary measures are supported on frozen
configurations. Indeed, each time a vertex updates, the number of disagreement
edges decreases strictly. Thus, if one starts from  a non-frozen shift-invariant
measure, then the density of disagreement edges decreases strictly over time,
which implies immediately that such a measure is not stationary (the argument can be made more precise along the lines of Remark~\ref{rem:Zd}).
\end{remark}

\subsection{The Ising Glauber dynamics on trees}

\label{sec:Ising}

Finally, we consider the Ising model on the directed tree
$\mathcal T=\mathcal T^3$.
The configurations of the model are again given by elements $\omega\in\{+,-\}^\mathbb V$,
where $\pm$ is identified with $\pm1$.
Consider some inverse temperature $\beta>0$ and some sequence 
of coupling constants
$(J_k)_{k\in\mathbb Z}\subset[0,\infty)$.
For this choice of parameters,
the \emph{edge energy} of the corresponding Ising model 
for an edge $e=xp(x)\in\mathbb E$ is given by
\[
   E(e)= -\beta J_k\omega_x\omega_{p(x)},    
\]
where $k$ is the unique integer such that $x$
belongs to the layer $\mathbb V_k$.
The Ising model is thus in principle defined on the undirected tree,
and the orientation of the edges appears only in the inhomogeneous choice of coupling constants.
The coupling constants $(J_k)_{k\geq 1}$ for positive layers do not matter,
but for the coupling constants of the negative layers $(J_k)_{k\leq 0}$
we require that
\begin{equation}
    \label{eq:J_requirement}
    J_{-k-1}-J_{-k}\geq k+o(k)
\end{equation}
for $k\geq 0$.
For simplicity, we assume in particular $J_{-k-1}> J_{-k}$
for all $k$,
which is automatically satisfied in any case for $k$ large.
For a concrete example, one may take
\begin{equation}
    \label{eq:choice}
    J_k:=\begin{cases}
        0&\text{if $k\geq 0$,}\\
        k^2&\text{if $k\leq 0$.}
    \end{cases}    
\end{equation}
    The edge energies induce a \emph{specification} (see~\cite{georgii2011gibbs} for a reference work on
the relevant technical constructions and which also considers the Ising model in detail).
   In simple words, to every
      finite subset $\Lambda$  of vertices and every $\bar\omega\in\{+,-\}^{\mathbb V}$ is associated a probability measure $\nu_{\Lambda,\bar\omega}$ on $\{+,-\}^{\mathbb V}$, where
      \[
      \nu_{\Lambda,\bar\omega}(\omega)\propto e^{-
        \sum_{e\sim\Lambda}E(e)}{\bf 1}_{\omega_{ \Lambda^c}=\bar\omega_{\Lambda^c}}\,
      \]
      the sum being extended to edges intersecting $\Lambda$ and $\omega_{\Lambda^c}$ denoting the restriction of $\omega$ to the complement of $\Lambda$.
A probability measure that is invariant under resampling of $\omega_\Lambda$ according to $\nu_{\Lambda,\omega}$ for \emph{all} finite $\Lambda$
is called a \emph{Gibbs measure}.

 We study the heat bath Glauber dynamics for the Ising model,
that is, the continuous-time Markov chain where the value of $\omega$
at each vertex $x$ is updated according to $\nu_{
  \{x\},\omega}$
whenever the Poisson clock at $x$ rings.

The times at which the clocks ring shall again be provided by the measure $\mathbb P$
introduced above.
This induces a natural coupling between this Markov chain and
a stationary state for the majority voter model,
which allows us to demonstrate that the stationary limit of the Markov chain
is not reversible and therefore not a Gibbs measure for the underlying Ising model.

\begin{theorem}
    \label{thm:ising}
    For $\beta>0$ sufficiently large, the Ising model on
    $\mathcal T^3$ defined above has a stationary measure for the
    Glauber dynamics which is not reversible and therefore not a Gibbs
    measure.
\end{theorem}

\begin{remark}
    The same result can be derived on trees $\mathcal T^d$
    with $d$ an odd integer above three.
\end{remark}

\begin{remark}
    Our primary model of interest in this study was the stochastic Ising model
    on trees.
    The majority voter model can be viewed as a zero-temperature limit
    of it: taking $\beta$ to infinity makes each vertex align
    with the majority of its children when it is updated, because $J_{-k-1}>J_{-k}$.
    The coalescing particles model, finally, provides an upper bound on how
    quickly new voter information flows along the graph in the voter model.
    More precisely, if the voters are not presented with new information,
    then the voters will certainly not change their opinion.
    Suppose that the
    coalescing particles model turns out not to have a nontrivial stationary state.
    This means essentially that the geometry of the information flow is such that eventually no
    new information flows through the model.
    This then implies that the voters are not presented new information,
    and therefore the voter model would have no non-frozen stationary state.
    Observe that the converse statement is not true:
    even if new information is occasionally presented, this does not guarantee the voters
    to change their opinion.
\end{remark}

\section{The coalescing particles model}

\label{sec:info-flow}

Assume the setting of Subsection~\ref{subsec:main:info}.
In this section we prove Theorem~\ref{thm:info}.
Let us first introduce boundary conditions
in the form of vertices that never change state.
For $x\in\mathbb V$, $V\subset\mathbb V$,
and $\omega\in\{0,1\}^\mathbb V$,
write 
\[
    \omega^{x,V}_y:=
    \begin{cases}
        \omega_y&\text{if $y\in V$,}\\
        \omega^x_y&\text{if $y\not\in V$.}
    \end{cases}
\]
Write $L^V$ for the generator corresponding to this local update,
and write $L^n$ for $L^V$ with $V=\cup_{k\leq n}\mathbb V_k$
whenever $n\in\mathbb Z$.
Unless indicated otherwise,
all Markov processes generated by $L$ and $L^n$ are always naturally coupled
with the sequences of Poisson clock times $(X_k^x)_{k\in\mathbb Z}$,
which means that the probability measure $\mathbb P$ contains all the randomness
for describing the evolution of the system.

Let us start with an observation on the Markov processes generated by $L$
and $L^n$.
First, each generator is monotone, that is it preserves the partial order of the system.
In other words, if $A^\omega=(A^\omega_t)_{t\geq 0}$ and $A^\eta=(A^\eta_t)_{t\geq 0}$ are Markov processes with the same generator and started from 
some configurations $\omega$ and $\eta$ respectively with $\omega\leq\eta$,
then $A^\omega_t\leq A^\eta_t$ for all $t$.

For any $s\in\mathbb R$, let $A^s=(A^s_t)_{t\in[s,\infty)}$ denote the Markov process
generated by $L$ using the coupling via the Poisson clocks, and started at time $s$ from the maximal configuration $\omega^1$ which has 
particles at all vertices.
If $s\leq s'$, then it is easy to see that
$A^s|_{[s',\infty)}\leq A^{s'}$,
because $A^s_{s'}\leq A^{s'}_{s'}=\omega^1$.
Define the process
$A^{-\infty}=(A^{-\infty}_t)_{t\in(-\infty,\infty)}$
by taking the pointwise limit
\[A^{-\infty}_t:=\lim_{s\to-\infty}A^s_t.\]
By definition, the law of $A^{-\infty}_t$ in $\mathbb P$ is
independent of $t$, which implies that it is stationary for the
generator $L$.  Moreover, the law of $A^{-\infty}_0$ in $\mathbb P$
must be the unique maximal stationary law for $L$, since it is obtained
as a limit of processes started from the maximal configuration
$\omega^1$.  It is also $\operatorname{Aut}(\mathcal T)$-invariant,
because the initial condition and the update rates
  are $\operatorname{Aut}(\mathcal T)$-invariant.  Of course, it is possible that the process $A^{-\infty}$ is
trivial, that is, that there is almost surely no particle at each
vertex.  We shall prove that this is not the case
(Theorem~\ref{thm:info}).

We now construct the process $A^{-\infty}$ as the limit of another sequence,
in a way which is more amenable to analysis.
Recall that the generator $L^n$ is identical to the generator $L$,
except that it does not update the state of vertices in $\cup_{k\leq n}\mathbb V_k$.
Introduce a new process $B^n$ which is a stationary Markov chain with generator $L^n$,
and which is identically equal to one on $\cup_{k\leq n}\mathbb V_k$.
Since this process has fixed boundary conditions on all layers $\mathbb V_k$ for $k$ 
sufficiently small,
it is easy to see through an argument involving coupling from the past~\cite{propp1997coupling} that 
this process is uniquely defined at all times---both positive and negative.
This is because particles flow only upwards from the base layers,
at which particles are always present.
Thus, if we fix a vertex and condition on all Poisson clocks
(for both positive and negative times),
then it is easy to see that its state at time zero does not depend 
on its state at time $-t$ when $t$ is larger than a random but almost surely finite time.
Write $B^n=(B_t^n)_{t\in\mathbb R}$ for this process.

Note that the family of processes $B^n$---indexed by $n$---live on the same probability space
(which defines the Poisson clocks),
and are therefore naturally coupled.
It is easy to see that in this coupling $B_t^n\leq B_t^m$
for any $n\leq m$ and at any time $t\in\mathbb R$.
In particular, the process $B^{-\infty}:=\lim_{n\to-\infty}B^n$ is well-defined as a pointwise limit.
Since each $B^n$ is stationary for $L^n$ and since $L^n$ converges pointwise to $L$,
this limit is stationary for $L$.

\begin{lemma}
    The processes $A^{-\infty}$ and $B^{-\infty}$ are almost surely equal.
\end{lemma}

\begin{proof}
    By taking into account the fact that $B^n$ is defined with all-one boundary conditions,
    it is straightforward to see that almost surely $B^n\geq A^{-\infty}$ for all $n$.
    Therefore almost surely $B^{-\infty}\geq A^{-\infty}$.
    But $B^{-\infty}$ is stationary for $L$ and $A^{-\infty}$
    is the maximal stationary process for the same generator,
    which implies that $B^{-\infty}\preceq A^{-\infty}$.
    This is possible only if $B^{-\infty}=A^{-\infty}$ almost surely.
\end{proof}

The fact that each process $B^n$ is stationary for $L^n$ makes the sequence 
$(B^n)_n$ easier to handle than the family $(A^s)_s$.
Although we cannot identify the law of $B^{-\infty}$ directly,
we can explicitly calculate some quantitative information.
For this purpose, let us first introduce the following family of events.
Recall that $(X^x_k)_{k\in\mathbb Z}$ denotes the sequence of times at which 
the Poisson clock at some vertex $x\in\mathbb V$ rings.
For any $n\in\mathbb Z\cup\{-\infty\}$, for any $s\leq t$, and for any $x\in\mathbb V$,
write $F_x^n(s,t)$ for the event
\[
    \left\{~\parbox{10.5em}{
        for $B^n$, there was a flow of particles towards $x$ between times $s$ and $t$
    }~\right\}
    =
    \left\{~\parbox{16em}{
        there is a $k\in\mathbb Z$ and a $y\in\mathcal C(x)$ such that $s\leq X^x_k< t$
        and such that $\lim_{u\uparrow X^x_k}B^n_u(y)=1$
    }~\right\}.
\]
Write $F_x(s,t)$ for $F_x^0(s,t)$ for brevity.
The objective is to calculate the probabilities of these events.
Note already that $\mathbb P(F_x(s,t))=\mathbb P(F_x(s+a,t+a))$ by stationarity of $B^0$,
and that $\mathbb P(F_{x_n}(s,t))$ is decreasing in $n$ due to monotonicity ({the path $(x_n)_{n\in\mathbb Z}$ was defined at the beginning of Section \ref{sec:def}}).
Define, for $n\geq 1$,
\[
    \rho_n:[0,\infty)\to[0,1],\,
    T\mapsto \mathbb P(F_{x_n}(0,T));
    \quad
    \bar\rho_n:=1-\rho_n.
\]

\begin{proposition}
    We have
    $\bar\rho_1(T)=e^{-T}$.
\end{proposition}

\begin{proof}
    For the process $B^0$, there are always particles at the base layer,
    and in particular at the vertex $x_0$.
    Thus, as soon as the Poisson clock at $x_1$ rings,
    there is a flow of particles from $x_0$ to $x_1$.
    The time to the first clock ring has the distribution of an exponential random variable.
\end{proof}

The key to calculating $\rho_n$ for all $n$ is provided by the following lemma.

\begin{lemma}
    \label{lemma:int_1}
    For all $n\geq 1$, we have
    \begin{equation}
        \label{eq:int_1}
        \rho_{n+1}(T)=
            \int_0^\infty e^{-s}ds
            \int_0^Te^{t-T}dt
            (1-\bar\rho_n(t+s)^2).
    \end{equation}
\end{lemma}

\begin{proof}
    Recall that $B^0$ is a stationary Markov chain defined for all times in $\mathbb R$.
    The formula is easy to derive from the probabilistic model.
    Let $-{\mathfrak s}$ and ${\mathfrak t}$ denote the last times before $0$ and $T$ respectively that the Poisson clock at the vertex $x_{n+1}$ rings.
    Thus, $-{\mathfrak s}$ is certainly negative (${\mathfrak s}$ is an exponential random variable of parameter $1$), and it is possible that $-{\mathfrak s}={\mathfrak t}$ if the clock at $x_{n+1}$ does not ring in the interval $[0,T]$.
    It is obvious that (up to sets of zero measure)
    \[
        F_{x_{n+1}}(0,T)=\{{\mathfrak t}\geq 0\}\cap(\cup_{y\in\mathcal C(x_{n+1})}F_y(-{\mathfrak s},{\mathfrak t})).
    \]
    Finally, for  fixed $a\leq b$, the event $\{{\mathfrak t}\geq 0\}$ and the events
    $(F_y(a,b))_{y\in\mathcal C(x_{n+1})}$ are all independent of one another, because particles flow only from children to parents and not vice versa.
    Therefore we have
    \begin{align*}
        \rho_{n+1}(T)
        &=
        \mathbb P(\{{\mathfrak t}\geq 0\}\cap(\cup_{y\in\mathcal C(x_{n+1})}F_y(-{\mathfrak s},{\mathfrak t})))
        \\
        &=
        \int_0^\infty e^{-s}ds
        \int_0^Te^{t-T}dt
        \mathbb P(\cup_{y\in\mathcal C(x_{n+1})}F_y(-s,t))
        \\
        &=
        \int_0^\infty e^{-s}ds
        \int_0^Te^{t-T}dt
        (
            1-\textstyle\prod_{y\in\mathcal C(x_{n+1})}
                (1-\mathbb P(F_y(-s,t)))
        )
        \\
        &=
        \int_0^\infty e^{-s}ds
        \int_0^Te^{t-T}dt
        (
            1-
                (1-\mathbb P(F_{x_n}(0,s+t)))^{|\mathcal C(x_{n+1})|}
        )
        \\
        &=
        \int_0^\infty e^{-s}ds
        \int_0^Te^{t-T}dt
        (1-\bar\rho_n(s+t)^2).
    \end{align*}
    In the second equality, we integrate out the distribution of ${\mathfrak s}$ and ${\mathfrak t}$ (using that $T-{\mathfrak t}$ is an exponential random variable of parameter $1$),
    conveniently absorbing the event $\{{\mathfrak t}\geq 0\}$ in the lower bound of the second integral.
    In all other equalities, we directly use the definitions of the quantities,
    as well as the independencies and symmetries of the events 
    under consideration.
\end{proof}

Define the endomorphism
\[
    \chi:\rho\mapsto
    \left(T\mapsto
    \int_0^\infty e^{-s}ds
    \int_0^Te^{t-T}dt
    (1-(1-\rho(t+s))^2)\right)
\]
on functions $\rho:[0,\infty)\to[0,1]$ which are increasing and upper bounded by $T\mapsto 1-e^{-T}$.
Clearly $\chi$ is monotone,
in the sense that $\chi(\rho)\leq \chi(\rho')$ if  $\rho\leq\rho'$ pointwise.
Since $\rho_1$ is the unique largest function in the domain of $\chi$,
it is immediate that $\rho_{n}=\chi^{n-1}(\rho_1)$ converges to the unique largest fixed 
point of $\chi$ as $n\to\infty$, which we call $\rho_\infty$.

\begin{lemma}
    \label{lemma:rho_infinity}
    For each $n\in\mathbb N\cup\{\infty\}$, we have 
    \(
        \rho_n(T)=\mathbb P(F^{-n}_{x_0}(0, T))
    \).
\end{lemma}

\begin{proof}
    For $n\in\mathbb N$,
    one can relate the law of $B^{-n}$ to that of $B^0$
    by applying the automorphism $\tau$ (defined at the beginning of Section \ref{sec:def}) $n$ times.
    This immediately yields the lemma.
    For the case $n=\infty$, observe that 
    $\rho_n$ decreases pointwise to $\rho_\infty$ 
    as $n\to\infty$.
    The result then follows by application 
    of the dominated convergence theorem, because the sequence of events $F_{x_0}^{-n}(0,T)$ is decreasing with intersection $F_{x_0}^{-\infty}(0,T)$.
    This in turn is a consequence from the fact that $B^{-\infty}$ is defined as the pointwise decreasing limit of $B^{-n}$.
\end{proof}

Thus, to prove Theorem~\ref{thm:info},
it suffices to identify $\rho_\infty$.
One can also wish to be more general and ask about all the fixed points of $\chi$.
Before doing so, we first rewrite the equations that a fixed point of $\chi$
must satisfy.

\begin{proposition}
    \label{pro:integral_transform}
    Consider the integral transform
    \[
        \rho\mapsto\left(T\mapsto
        \int_0^\infty e^{-s}ds
        \int_0^Te^{t-T}dt
        f(\rho(s+t))\right),
    \]
    for functions $\rho:[0,\infty)\to\mathbb R$
    where $f$ is some polynomial.
    Then all fixed points of this integral transform
    are twice  differentiable,
    and satisfy
    \(
        \rho''=\rho-f(\rho)
    \).
\end{proposition}

\begin{proof}
    First differentiate both sides of the equation to make the integral over $t$ disappear.
    Then substitute $s+T$ by $s$ in the resulting expression and differentiate again,
    to obtain the desired equation.
\end{proof}

This is convenient, because we may interpret the equation $\rho''=\rho-f(\rho)$
as describing the trajectory of a one-dimensional Newtonian particle moving in some potential $V$.
In particular, we write
\[
    \rho''=-V'(\rho);\quad V'(\rho)=f(\rho)-\rho.
\]
Substituting the polynomial for $f(p)=(1-(1-\rho)^2)$ and taking the primitive yields
\(
    V(\rho)=\frac12\rho^2-\frac13\rho^3.
\)
This is good news, because
$V$ is non-decreasing on $[0,1]$,
and because $V'$ vanishes at $0$ and $1$.
Thus, this tells us that there are exactly two solutions
$\rho:[0,\infty)\to[0,1]$ which satisfy $\rho(0)=0$;
one that is identically zero,
and one which tends to one as $T$ tends to infinity.
It can be shown that $\rho_\infty$, the latter solution, equals
\[
    \rho_\infty(T)=1+6\frac{\beta e^{-T}}{(\beta e^{-T}-1)^2}
    ;\quad\beta=\sqrt 3 - 2,
\]
by plugging this solution into the differential equation $\rho''=\rho-f(\rho)$.
This proves Theorem~\ref{thm:info}.

If the interest is in the same model but on $\mathcal T^d$ rather than $\mathcal T^2$,
then one may replace the power $2$ in the equations by $d$.
This leads to another ordinary differential equation which has $\rho_\infty$
as its nontrivial solution.

\section{The majority voter model}

To prove Theorem~\ref{thm:nonfrozenvoters} we proceed as in the previous subsection by introducing a process $B$ 
with boundary conditions.
Boundary conditions, however, are enforced in a slightly more convolved way.
The process $B$ is a Markov chain on the state space $\{+,-\}^\mathbb V$
with the following update rules.
If a clock rings at a vertex $x$ in $\cup_{k>0}\mathbb V_k$,
then the current configuration $\omega$ is replaced by $\omega^x$.
However, if a clock rings at some vertex $x$ in $\cup_{k\leq 0}\mathbb V_k$,
then the new opinion of $x$ is instead determined by flipping 
a fair, independent $\pm$-valued coin.
{
This means that those vertices act as a sort of boundary whose random behaviour is imposed,
and this boundary acts as an independent source of randomness.}
We shall make a slight abuse of notation by pretending that $B$ lives in the measurable space corresponding 
to $\mathbb P$,
even though the process formally requires extra auxiliary randomness
to decide on the outcomes of the coin flips.
As in the previous section,
the process $B=(B_t)_{t\in\mathbb R}$ can be defined at all times 
through standard arguments involving coupling from the past~\cite{propp1997coupling}.
For this argument it is essential to observe that the state of a vertex depends on the past of only
finitely many other vertices (namely itself and its descendants which are not in the interior of the boundary).
Thus, the initial distribution $B_0$ is stationary for the dynamics described above.

Recall the definition of the automorphism $\tau$ from Section \ref{sec:def} and denote by $\mathbb P_k$ (with corresponding expectation $\mathbb E_k$) the law of $\tau^k(B)$ as a measure over trajectories, where the stationary process $\tau^k(B)$ is the process $B$ after $k$-fold application of $\tau$. Let $L^k$ denote the generator of $\tau^k(B)$ and $\mu_k$ its stationary distribution.
    Note that each such process is coupled to the Poisson clocks of $\mathbb P$ in a natural way.
    The state space $\{+,-\}$ is compact and we chose the product topology on $\{+,-\}^{\mathbb V}$, and therefore  we may always extract subsequential limits from sequences of (linear combinations of) such measures $\mu_k$.
   \begin{lemma}\label{lemma:430}
Any subsequential limit of the Ces\`aro sum
\begin{equation}
    \label{eq:cesaro}
    \frac1n\sum_{k=n}^{2n-1}\mu_k
\end{equation}
is an $\operatorname{Aut}(\mathcal T)$-invariant stationary probability measure for the Markov process with generator $L$.
   \end{lemma}
   \begin{proof}
     Call $\mu$ one such limit point. To check stationarity, we need that $\mu(L f)=0$ for every fixed, local (and therefore bounded) function $f$. We know that $\mu_k(L^k f)=0$. On the other hand, $L^k f=L f$ whenever $f$ is measurable with respect to the spins on $\cup_{j>-k+1}\mathbb V_j$, which holds for $n\le k\le 2n$, if $n$ is large enough (depending on the function $f$).
     Notice that this follows immediately from the definition of the two generators.
     To see that $\mu$ is $\operatorname{Aut}(\mathcal T)$-invariant, note first that any automorphism $\phi$ of $\mathcal T$ can be written as the composition of $\tau^z$ for some $z\in\mathbb Z$ and of some automorphism $\psi$ that maps each level $\mathbb V_j,j\in\mathbb Z$ into itself. Secondly, by construction of the process $B$, $\mu_k$ is invariant under such automorphism $\psi$. Therefore (with the convention that for an event $A$, $\phi(A)$ denotes the event $\{\omega:\phi^{-1}(\omega)\in A\}$, and
where the limit is meant to be taken along the subsequence that defines $\mu$)
     \begin{multline}
       \mu(\phi(A))=\lim_{n\to\infty}\frac1n\sum_{k=n}^{2n-1}\mu_k(\psi\circ\tau^z(A))=\lim_{n\to\infty}\frac1n\sum_{k=n}^{2n-1}\mu_k(\tau^z(A))\\=\lim_{n\to\infty}\frac1n\sum_{k=n-z}^{2n-1-z}\mu_k(A)=\mu(A)\,.
     \end{multline}
   \end{proof}
Our strategy will simply be to prove that even far away from the base layer $\mathbb V_0$,
the opinions of the vertices change sufficiently frequently,
so that this subsequential limit is non-frozen.
Remark that for the model of coalescing particles,
the monotonicity property allowed us to take limits rather than subsequential limits.
Since we lack a useful form of monotonicity in the current setting,
it remains open to ask if the subsequential limit obtained above 
is the ``unique maximal state'' for an appropriate partial ordering.

In analogy to the method in the previous section,
we define, for each $n\geq 0$,
\[
\bar\rho_n:[0,\infty)\to[0,1],\,T\mapsto\mathbb E[B_0(x_n)B_T(x_n)]
  =\mathbb E_n[B_0(x_0)B_T(x_0)],
  \]
where for the sake of calculating the expectation,
the values $\pm$ are identified with $\pm1$.
We also set $\rho_n:=1-\bar\rho_n$.
Observe  that
\[
    \bar\rho_0(T)=\mathbb P(\text{the clock at $x_0$ does not ring in $[0,T]$})=e^{-T}.    
    \]
    The first equality holds because the spin at $x_0$ stays the same up to the first positive time that the Poisson clock at $x_0$ rings, while at the first clock ring the spin is resampled as an independent centred random variable.
    The probability that the clock does not ring in the interval $[0,T]$ is exactly $e^{-T}$.
We will calculate $\rho_n$ for any $n>0$ by deriving an integral equation 
similar to~\eqref{eq:int_1}.

First introduce an auxiliary function.
Suppose we flip three fair, independent coins,
and write $x$ for the majority outcome.
Then, each coin is turned over with probability $p$,
independently of anything else (by \emph{turning over} we mean that we change the outcome of the coin;
not that we resample its value).
Write $y$ for the majority outcome after this procedure,
and write $M(p)$ for the probability that $x\neq y$.

\begin{proposition}
    We have $M(p)=\frac32p-\frac32p^2+p^3$, which is increasing in $p$.
\end{proposition}

\begin{proof}
It is straightforward to work out that
\begin{enumerate}
    \item $M$ must be a cubic polynomial in $p$,
    \item $M(0)=0$, $M(\frac12)=\frac12$, and $M(1)=1$,
    \item $M'(0)=\frac32$,
\end{enumerate}
which uniquely identifies $M$ as being the asserted polynomial.
{To see that $M(1/2)=1/2$, we remark that when $p=1/2$, turning over a coin with probability $p$ has the same effect (in distribution) as resampling it independently (for $p\ne 1/2$, these two operations give a different result). Therefore, all three coins are independently resampled, and the new majority is independent from the original one.}
To see that $M'(0)=\frac32$,
observe that $M'(0)$ equals the number of coins times 
the probability that turning a single coin changes the majority vote.
The latter probability equals $\frac12$ in the case of three coins.
\end{proof}

\begin{lemma}
    For any $n\geq 0$, we have
    \begin{equation}
        \label{eq:int_2}
        \rho_{n+1}(T)=
        \int_0^\infty e^{-s}ds
        \int_0^Te^{t-T}dt
        \,2M({\textstyle\frac12}\rho_n(t+s)).    
    \end{equation}
\end{lemma}

\begin{proof}
    The proof is---in broad strokes---the same as the proof of 
    Lemma~\ref{lemma:int_1}.
    The integrals appear in the same way;
    it is only the change in integrand that we must justify here.
    Note first that
    \[
        \rho_{n+1}(T)=2\mathbb P(B_0(x_{n+1})\neq B_T(x_{n+1})).
    \]
    If we write $-\mathfrak{s}$ and $\mathfrak{t}$ for the same times as in Lemma~\ref{lemma:int_1},
    then we are in fact comparing the outcomes of the majority votes 
    of the children of the vertex $x_{n+1}$ at times $-\mathfrak{s}$ and $\mathfrak{t}$.

    Note that at time $-\mathfrak{s}$, the restriction $B_{-\mathfrak{s}}|_{\mathbb V_n}$
    of $B_{-\mathfrak{s}}$ to the single layer $\mathbb V_n$ has the distribution of independent 
    fair coin flips.
    Between time $-\mathfrak{s}$ and time $\mathfrak{t}$, each coin has been turned 
    with probability $\frac12\rho_n(\mathfrak{t}+\mathfrak{s})$,
    independently of anything else.
    This justifies the new choice of integrand 
    $2M(\frac12\rho_n(\mathfrak{t}+\mathfrak{s}))$ in~\eqref{eq:int_2}.
\end{proof}

Define the endomorphism
\[
    \chi:\rho\mapsto
    \int_0^\infty e^{-s}ds
    \int_0^Te^{t-T}dt
    \,2M({\textstyle\frac12}\rho(t+s)),
\]
on functions $\rho:[0,\infty)\to[0,1]$ which are increasing and
bounded above by $T\mapsto 1-e^{-T}$.  Clearly $\chi$ is monotone,
since $M$ is also a monotone function.  Since $\rho_0$ is the unique
largest function in the domain of $\chi$, it is immediate that
$\rho_{n}=\chi^{n}(\rho_0)$ converges to the unique largest fixed
point of $\chi$ as $n\to\infty$, which we call $\rho_\infty$.
Proposition~\ref{pro:integral_transform} asserts that the fixed points
of $\chi$ may be found by solving
\[
    \textstyle
\rho''
=
\rho-f(\rho)
=
\rho-2M(\frac12\rho).
\]
We interpret this equation again as $\rho$ being the position of  a one-dimensional particle moving in a potential,
and write
\[  
    \textstyle
    \rho''=-V'(\rho); \qquad
    V(\rho)=\frac14\rho^2-\frac14\rho^3+\frac1{16}\rho^4.
\]
Again, we observe that $V$ is non-decreasing on $[0,1]$,
with $V'$ vanishing at $0$ and $1$.
This means that there are two distinct solutions,
one being $\rho\equiv 0$, and the other one being the solution $\rho_\infty$
introduced before.
Even without knowing the closed-form formula for $\rho_\infty$,
we can see that $\bar\rho_\infty$ decays exponentially fast in time,
since $-V'(\rho)$ behaves as $\frac 14(1-\rho)$ as $\rho\uparrow1$.

\begin{proof}[Proof of Theorem~\ref{thm:nonfrozenvoters}]
    Let $\mu$ denote any subsequential limit of the Ces\`aro
    sum~\eqref{eq:cesaro} as $n\to\infty$ in the topology of local
    convergence and $B$ denote the stationary and automorphism-invariant (thanks to Lemma \ref{lemma:430}) Markov process with
    generator $L$, initial measure $\mu$, and law $\mathbb P$. 
    Then,
\[
    \mathbb E[B_0(x)B_T(x)]=\lim_{n\to\infty} \mathbb E_n[B_0(x)B_T(x)]=\bar\rho_\infty(T)
\]
for any vertex $x$ and for any $T\geq 0$, where the limit is taken along the subsequence.
It was proved above that $\bar\rho_\infty$ decays exponentially fast,
which means that the opinion of each vertex mixes exponentially fast in time.
The latter implies also that $B$ is not supported on frozen (i.e., absorbing) states.
\end{proof}

\section{The Ising model}

\begin{proof}[Proof of Theorem~\ref{thm:ising}]
    In this proof we assume the choice~\eqref{eq:choice} of coupling constants for simplicity;
    the general case follows by the same argument.
Now construct two processes $(B_t)_{t\geq 0}$ and $(F_t)_{t\geq 0}$
taking values on $\{+,-\}^\mathbb V$ in the same probability measure $\mathbb P$
as follows.
\begin{itemize}
    \item The two processes start from the same configuration, which is sampled
    from the stationary distribution of the voter model
    constructed in the proof of Theorem~\ref{thm:nonfrozenvoters}.
    \item The probability measure $\mathbb P$ features a number of independent Poisson clocks
    at the vertices of the graph.
    \item Whenever a clock at some vertex $x$ rings at time $t$, then the value of $B_t(x)$ is replaced
    by the outcome of a majority vote of its children.
    Notice that this implies that $(B_t)_{t\geq 0}$ is a stationary Markov process with the generator $L$ introduced in Subsection~\ref{sec:voter}.
    \item Whenever a clock at some vertex $x$ rings at time $t$, then the value of $F_t(x)$ is replaced
    by resampling the value of the spin at $x$ according to the kernel $\nu_{\{x\},F_t}$ of the Ising model specification corresponding to $\{x\}$.
    Notice that more auxiliary randomness is required to define the outcome of this resampling operation.
\end{itemize}
The idea is that the coupling constants in the Ising model are chosen such that
the majority vote and the Ising update at some vertex $x$ agree with high likelihood, whenever
the two configurations are already equal at the children of $x$.
We shall formalise this by introducing a new process that contains the set of vertices
at which $B$ and $F$  disagree, then prove some simple upper bounds on the size of this new process.

Define $(D_t)_{t\geq 0}$ by $D_t:=1_{\{B_t\neq F_t\}}\in\{0,1\}^\mathbb V$.
We think of $D$ as a disagreement percolation process,
and we also view $D_t$ as a random subset of $\mathbb V$.
The goal of this section is to prove that $D$ does not grow very large in a precise sense,
so that $B$ and $F$ remain perfectly coupled on large portions of the graph.

First define
\[
    D_t':=D_t\cup \{x\in\mathbb V:\text{$D_t$ contains some descendant of $x$}\}.
\]
Suppose that the children of $x\in\mathbb V_k$ do not belong to the random set $D'$ at some time $t$,
where $k\leq 0$.
In particular, this means that $B$ and $F$ agree at the children of $x$.
The probability that a Glauber update would disagree with a majority vote at the children
is upper bounded by
\begin{equation}
    \label{eq:upperbounds}
    \frac
        {e^{-2\beta (2|k|+1)}}
        {1+e^{-2\beta (2|k|+1)}}
    \leq
    e^{-2\beta (2|k|+1)}
    \leq \frac12,
\end{equation}
where for the second inequality we impose that $\beta\geq 1$. In fact, the most unfavorable case is that where two children have a different sign from the third, and we recall that $J_{-k-1}-J_{-k}=2|k|+1$.
Using these inequalities, it is easy to define a new process $D''$ which stochastically dominates
$D'$.
The process $(D''_t)_{t\geq 0}$ is defined as follows.
\begin{itemize}
    \item At all times $t\geq 0$, we have $\cup_{k>0}\mathbb V_k\subset D''_t$, with equality at time $t=0$.
    \item If a Poisson clock at a vertex $x\in\cup_{k\leq 0}\mathbb V_k$ rings at time $t$,
    then this vertex and all its ancestors are added to $D''$ with probability $e^{-2\beta (2|k|+1)}$.
    Notice that vertices are added to $D''$ with a higher rate than to $D'$ due to the inequalities
    in~\eqref{eq:upperbounds}.
    \item If a Poisson clock at a vertex $x\in\cup_{k\leq 0}\mathbb V_k$ rings at time $t$
    and if all children of $x$ are not in $D''$,
    then $x$ is removed from $D''$ with probability $1/2$.
    Notice that vertices are removed from $D''$ with a lower rate than from $D'$
    due to the inequalities in~\eqref{eq:upperbounds}.
    \item The two operations above are coupled such that they never both occur at the same time,
    which is possible because the sum of the two probabilities does not exceed one.
    This implies that if a Poisson clock at a vertex $x\in\cup_{k\leq 0}\mathbb V_k$ rings at time $t$
    and if all children of $x$ are not in $D''$,
    then with probability $1/2-e^{-2\beta (2|k|+1)}$, nothing is done.
\end{itemize}
Notice that $D''$ is also closed under adding ancestors.
By the dominations described above, it is easy to see that $D''_t$ stochastically
dominates $D'$,
and in fact it is easy to construct all processes on the same probability space
such that $D_t\subset D'_t\subset D''_t$ at all times almost surely (by adding more auxiliary randomness if necessary).

Now let $N_t^k$ denote the random number of descendants of $x_k$
that belong to $D''_t$,
where $k\leq 0$.
Notice that $N^k$ is stochastically dominated by the Markov process on $\mathbb Z_{\geq 0}$
which:
\begin{enumerate}
    \item Jumps up by $i\in\mathbb Z_{> 0}$ with rate $3^i\cdot e^{-2\beta (2(|k|+i)+1)}$,
    \item Jumps down by $1$ with rate $1/2$ as long as its value is positive.
\end{enumerate}
Recall that $\beta\geq 1$, which is important because it implies that the first rate is 
exponentially decreasing with $i$.
It is a straightforward exercise in the theory of Markov chains that,
for $|k|$ large enough,
this process has a unique stationary probability distribution with exponential tails.

This means that for $k\in\mathbb Z_{\leq 0}$, the quantity
\begin{equation}
    \label{eq:uniform_bound}
    \sup_{t\geq 0}\mathbb P(x_k\in D''_t)
\end{equation}
is exponentially decreasing in $|k|$.
Since $D\subset D'\subset D''$, this also implies the bounds on $D$ and $D'$ that we sought after.

Write $\mu_t$ for the joint distribution of $(B_t,F_t)$, observing that
$D'_t$ depends deterministically on this pair.
By compactness,
the sequence of probability measures
\[
    \nu_T:=\frac1T\int_T^{2T}dt \mu_t
\]
has some subsequential limit $\nu$ as $T\to\infty$.
Write $(\tilde B,\tilde F)$ for the random objects with initial distribution given by the probability measure $\nu$,
and write $\tilde D$ and $\tilde D'$ for the derived percolations.
This means that $\tilde B$ and $\tilde F$ are stationary distributions for the voter
dynamics and the Ising model dynamics respectively, coupled via the limit coupling just described.
The percolation $\tilde D$ is the disagreement percolation of those two processes,
and $\tilde D'$ contains $\tilde D$ and ancestors of vertices in $\tilde D$.
The equation defining $\nu_T$ should be thought of as a continuous version of the
Ces\`aro sum before,
and it is straightforward to verify as in the proof of Lemma \ref{lemma:430} that indeed $\nu$
is stationary for the joint dynamics on the pair $(B_t,F_t)$ defined above.
Of course, the uniform bound in~\eqref{eq:uniform_bound}
passes to the limit, implying that
\[
    \nu(x_k\in \tilde D')   
\]
decays exponentially in $|k|$ over $k\in\mathbb Z_{\leq 0}$.
But if we now start our joint dynamics on the stationary pair $(\tilde B,\tilde F)$,
then this exponential decay implies that at a fixed time the two processes are equal on large
portions of the space (by summability of $\sum_{k=0}^\infty a^k$ for $|a|<1$),
which in turn means that the stationary distribution 
for the Ising model Glauber dynamics inherits its 
non-reversibility from the voter model.
\end{proof}

\section*{Acknowledgements}
The authors thank Aernout van Enter and Elchanan Mossel for pointing out several
relevant articles in the literature,
and the anonymous referees for numerous useful suggestions for improvement.
P.\ L.\ thanks F.\ T.\ and the TU Wien for their hospitality during
several visits.
F.\ T.\ gratefully acknowledges financial support of the Austria Science
Fund (FWF), Project Number P 35428-N.
This project has received
funding from the European Research Council (ERC) under the European Union's Horizon
2020 research and innovation programme (grant agreement No. 757296).

\bibliographystyle{amsalpha}
\bibliography{references.bib}

\end{document}